\documentclass[10pt]{amsart}

\usepackage{graphicx}
\usepackage{amssymb}
\usepackage{amsmath,amsthm,amscd,amssymb}
\usepackage{latexsym}

\usepackage{amsfonts}
\usepackage{amsthm}
\usepackage{enumerate}
\usepackage{epsfig}

\numberwithin{equation}{section}

\theoremstyle{plain}
\newtheorem{theorem}{Theorem}[section]
\newtheorem{lemma}[theorem]{Lemma}
\newtheorem{corollary}[theorem]{Corollary}

\newtheorem{remark}[theorem]{Remark}

\newtheorem{case[theorem]}{Case}

\title{On a class of curved flag multipliers}
\sf{\author{Hadi Jor\'{a}ti}}

\newcommand{\R}{\mathbb{R}}
\newcommand{\integral}{\int_{-\infty}^{\infty}}

\begin{document}

\maketitle

\begin{abstract}

We give a Mikhlin multplier theorem for a class of nonhomegenous dilations in plane,
and study singular kernels 
that have a flag type singularity along the parabola $y=cx^2$.
We show that the multiplier of these operators consists of two distinct parts,
one which is essentially the Fourier transform of the flat kernel,
and another, which is characterized by its highly oscillatory factor
that signifies the curvature of the singular support of the kernel.
This characterization can be used to study the behaviour
of this class on various function spaces. For example
the $L^p$ boundedness of the corresponding operator on
$1<p<\infty$ will be a trivial corollary..

\end{abstract}


\section{Introduction}

The standard theory of singular integral operators has been generalized to apply to certain wider
classes of operators in two main directions.  
One, in connection to the singular Radon transform operators, whose kernels are supported on
curved subvarieties. The pinnacle of this study was the paper of Christ, Nagel, Stein, Wainger \cite{annals}.
The other generalization, arises in connection to the product theory of singular intgeral operators,
where the singularity is carried on a flag of linear subvarieties, as suggested by speculations on problems 
in several complex variables. The paper of Nagel, Ricci, Stein \cite{flag} offers an extensive study.

In this article, we combine the two directions, in a special case, to study a class of operators that 
have singularities along a specific curved flag. These operators can also be viewed as a generalization 
of the classical operator of the Hilbert transform along the parabola, in the context of curved flag kernels.
We study the Fourier transform of such kernels, obtaining the exact form in one region, characterized by 
a very specific oscillatory factor, and establishing the smoothness in the rest of the plane.
As a corollary of the form of the multiplier, we can also show the $L^p$ boundedness
of the corresponding class of operators for $1<p<\infty$. The paper of Secco \cite{secco}
contains a direct proof of this corollary.


\bigskip
We will use the theory developed by Nagel, Ricci, and Stein cf. \cite{flag}, 
in the special case of $\mathbb{R}^2$
where we have assigned the homogeneous dimensions 1,2 to $x$ and $y$ directions respectively.
Also, in this special case we will be dealing with the following two filtrations of $\mathbb{R}^2$:
$$
\textrm{flag } \mathcal{F}_1: \qquad 0=V_0 \subset V_1=\mathbb{R}\times{0} \subset V_2=\mathbb{R}^2,
\qquad
\textrm{flag } \mathcal{F}_2: \qquad 0=V_0 \subset V_1={0}\times \mathbb{R} \subset V_2=\mathbb{R}^2
$$

So in this setting a {\em product kernel} is a distribution $M$ on $\mathbb{R}^2$ that coincides with a $C^{\infty}$ function $M$
away from the coordinate axes, which also satisfies the following decay and cancellation condirions:

(1)(Differential inequalities) 

For each pair of integer indices $\alpha, \beta$, there is a constant $C_{\alpha, \beta}$
such that
$$
| \partial_{x}^{\alpha}\partial_{y}^{\beta} M(x,y)| \leq C_{\alpha,\beta} |x|^{-1-\alpha} |y|^{-1-\beta}$$

(2)(Cancellation conditions) for any normalized bump functions $\phi, \psi$
each of the two families $\int M(x,y)\phi(Rx)dx$ and $\int M(x,y)\psi(Ry)dy$
are uniformly bounded families of CZ kernels on $\mathbb{R}$

A flat {\em flag kernel} with respect to flag $\mathcal{F}_1$ is a distribution $M$ on $\mathbb{R}^2$
that coincides with a $C^{\infty}$ function $M$ away from the $x$-axis, which satisfies
the following differential inequalities: for $y \neq 0$,
$$
| \partial_{x}^{\alpha}\partial_{y}^{\beta} M(x,y)| \leq C_{\alpha,\beta} (|x|+|y|^{1/2})^{-1-\alpha} |y|^{-1-\beta}$$

as well as the same cancellation conditions we imposed on product kernels. \cite{flag}

So, basically, a flag kernel with respect to this flag coincides with a product kernel 
in the region where $|x|^2 \geq c|y|$, and a Calderon-Zygmund kernel with respect to 
the corresponding weights in the remainder of the plane: $|y| \geq c|x|^2$.

Notice that any product kernel can be written as a finite sum of flag kernels (cf. \cite{flag})
and that the sum of a flag kernel with another flag kernel with respect to a coarser 
filtration is a flag kernel with respect to the finer flag. In particular,
the sum of a Calderon-Zygmund kernel and any flag kernel is a flag kernel of the same kind.


\bigskip

Our main results are the following two dual theorems that characterize a class of singular kernels
that behave naturally under a specific class of nonhomegenuous dilations in $\mathbb{R}^2$, 
giving a new kind of Mikhlin multiplier type condition.
There is also a rather precise asymptotic relation between the different terms.

\medskip
\begin{theorem}
If $M$ is a flat flag kernel with respect to the flag $\mathcal{F}_1$
and $K$ is the curved version, $K(x,y)=M(x,y-{c_0}x^2)$, and $m$ the Fourier transform of $K$, 
then $m$ is a bounded function which is smooth away from $\xi$-axis, 
and there is a natural decomposition for it:


\begin{equation}
m(\xi,\eta)={\tt F.T. } \{M(x,y-{c_0}x^2)\}(\xi,\eta)=
L_{1}(\xi,\eta)+\Phi(\xi,\eta)e^{ic'\frac{\xi^2}{\eta}}\frac{\eta^{1/2}}{\xi}L_{2}(\xi,\eta)
\end{equation}

where each $L_{i}$ is a flat flag multiplier with respect to the flag $\mathcal{F}_1$,
$\Phi$ is a smooth cut off function supported on the region where $|\xi| \geq c |\eta|^{\frac{1}{2}}$.
\end{theorem}

\medskip
\begin{theorem}
If $\ell$ is a flag multiplier with respect to flag $\mathcal{F}_1$, then 
for the corresponding oscillatory expression we have:
\begin{equation}
{\tt Inverse  F.T. }
\{\Phi(\xi,\eta)e^{ic'\frac{\xi^2}{\eta}}\frac{\eta^{1/2}}{\xi}\ell(\xi,\eta)\}(x,y)
=M_2(x,y) + M_1(x,y-{c_0}x^2)
\end{equation}

where each $M_i$ is a flag kernel 
with respect to flag $\mathcal{F}_i$.
\end{theorem}

\medskip
There is an asymptotic relation between the different parts in both directions. 

Introducing parabolic polar co\"ordinates $(a,\delta)$  as 
$a=\frac{\xi}{\eta^{1/2}}, \quad \delta=\eta^{1/2}$ for $\xi, \eta > 0$,
we write for two functions $m$ and $\ell$ on the $\xi,\eta$ plane 
$m \sim \ell$ as $a \to \infty$, 
if we have, for all positive integers $r$:
$$
\partial^{r}_{a}  (m-\ell)(\xi,\eta)
=O\left( (1+|a|)^{-r-1}\right)
\qquad
\textrm{as}
\quad
a \to \infty
$$
with constants not depending on $\delta$.

\bigskip
\begin{remark}(Asymptotic behaviour)
There is the following asymptotic relations between $L_1, L_2$ and $M_1, M_2$
for $\frac{\xi^2}{\eta} \to \infty$:

\begin{equation}
L_1(\xi,\eta) \sim {\tt F.T.} \{ \frac{1}{x}L_2(\eta x, \eta)\}(\xi)
\hspace{2cm}
\widehat{M_1}(\xi,\eta) \sim  \widehat{M_2}(\xi,\eta) 
\end{equation}

\end{remark}

It is worth noting that one can't hope for anything better than the asymptotic relation since the
exact form will depend on many factors such as the specific cut off function
$\Phi$ and therefore cannot be contorlled.


\vspace{6mm}

\section{general setting and basic lemmas}

Throughout these pages $c$ 
means a constant whose value is independent of all the variables and different appearances 
may refer to different constants.

We call a $C^{\infty}$ function $\phi$ on $\R^n$ a {\it normalized bump function} 
if it is supported in the unit ball and all its $C^r$ norms, for $0 \leq r \leq 1$, are bounded by $1$.


A {\it Calderon-Zygmund} distribution on $\mathbb{R}$, is a distribution $k$ that equals a function $k$ away from the origin
which satisfies the following decay and cancellation conditions:

i) (decay) \hspace{3cm}
$|\partial_{x}^{\alpha}k(x)| \leq C_{\alpha} |x|^{-1-\alpha}$

ii) (cancellation) 
$|\int k^{\delta}(x)\phi(x)dx| \leq C$
with $C$ not depending on $\delta$ or $\phi$ (normalized bump function)
and
$k^{\delta}(x):=\delta^{-1}k(\delta^{-1}x)$ represents dilations of $k$.
The best such constants $C, \{C_{\alpha}\}$ are called {\it CZ seminorms}.


We call a family of CZ kernels $\{k_i\}_i$, a uniformly bounded family, 
if for any finite number of CZ seminorms $\|.\|_{\alpha}$, the family 
$\{\|k_i\|_{\alpha} \}_{\alpha \leq N}$
is uniformly bounded:
$
\|k_i\|_{\alpha} \leq C$  for 
$\alpha \leq N$, with $C$ independent of $i$.




It is easy to check that
if $k$ is a Calderon-Zygmund kernel on $\R$, then, for each $\alpha$, the family 
$\left\{ \delta^{\alpha} \partial_{\delta}^{\alpha} k^{\delta} \right\}_{\delta>0}$
is a uniformly bounded family of CZ kernels.

%
%
%
%
%
%
%
%
%
%

Also, the following trivial corollary  will be handy:
If $m$ is a Calderon-Zygmund multiplier on $\R$ then, for each $\alpha$, the family 
$\left\{ \delta^{\alpha} \partial_{\delta}^{\alpha} m_{\delta} \right\}_{\delta}$
is a uniformly bounded family of Calderon-Zygmund multipliers. 


\bigskip
As a model case we first
study a generalization of the operator of Hilbert transform on the parabola 
where we have replaced the homogeneous kernel 
$\frac{1}{t}$
with a general $1$-dimensional CZ kernel $k$:
$$
T_k(f)(x,y):= f*K_{k}(x,y)$$

where $K_{k}$ is defined as: $K_{k}(\phi):=k(\phi(x,x^2))$.
So $\widehat{K_{k}}(\phi)=k(\widehat{\phi}(\xi,\xi^2))$.\\
But, every Calderon-Zygmund kernel in dimension one
is in the following form (see e.g. \cite{bigstein}):
$k=\lim_{\epsilon \to 0, N \to \infty} k_{\epsilon,N} +c\delta$
where 
$\delta$ is the $\delta$-function at zero, and $k_{\epsilon,N}$ are $C^{\infty}$ functions 
supported in the shell between radii $\epsilon$ and $N$, 
with uniformly bounded CZ semi-norms.


Hence
$$
\widehat{K}(\phi)= \lim_{\epsilon \to 0, N \to \infty} 
\int \left\{ \int\!\!\int e^{-2{\pi}i(tx+{t^2}y)}\phi(x,y) dxdy \right\} 
k_{\epsilon,N}(t) dt + c\widehat{\delta}(\phi)$$

But 
$$
\widehat{\delta}(\phi)=\delta(\widehat{\phi}(\xi,\xi^2))=\widehat{\phi}(0,0)
=\int\!\!\int \left\{ \int e^{-2{\pi}i(tx+{t^2}y)} \delta(t) dt \right\} \phi(x,y) dxdy
=\int\!\!\int \phi(x,y) dxdy$$

Thus the whole expression representing the Fourier transform of $K_{k}$ is equal to:
$$
\int\!\!\int \left\{ \int e^{-2{\pi}i(tx+{t^2}y)} 
\{ \lim_{\epsilon \to 0, N \to \infty}  k_{\epsilon,N}(t)+c\delta \} dt 
\right\} \phi(x,y)dxdy$$

So if we call the multiplier of this operator $m_k(\xi,\eta)$, 
in slack notation, $m_k$ is given by the following integral: 
$$
m_k(\xi,\eta)=\int_{-\infty}^{\infty} e^{-2{\pi}i(t.\xi+{t^2}.\eta)}k(t)dt$$

By this we mean, 
$m_{k}(\xi,\eta)=\lim_{\epsilon \to 0, N \to \infty} m_{\epsilon,N}+c$,
where $m_{\epsilon,N}$ is the Fourier transform of $k_{\epsilon,N}$.


\begin{lemma}

If $P$ is a polynomial of degree $d$, and $\mathcal{K}:=\{k_\alpha\}_{\alpha}$ is a uniformly bounded family of CZ kernels,
\begin{displaymath}
|\int e^{iP(t)} k_{\alpha,\epsilon,N}(t)dt| \leq c_{d, \mathcal{K}} 
\end{displaymath}

where $c_{d, \mathcal{K}}$ is a constant depending only on $d$, the degree of $P$, and the uniform bounds
of the CZ seminorms of $\{k_\alpha\}$ but not on the coefficients of $P$, or
$\epsilon, N$. 
\end{lemma}

\noindent
An earlier version is due to Stein and Wainger \cite{studia-70}.


\noindent
{\bf Proof.}
We can see, without loss of generality, that we can assume $P$ is monic. Indeed, if it is not,
we can perform a change of variables, rescaling the variable appropriately, to make it monic, while
this change of variables replaces any kernel $k_{\alpha,\epsilon,N}$ with a dilated version of it, 
and hence leaving the CZ seminorms unchanged.

Break up the integral into two parts: the part near zero, and the part near infinity.
First we treat the part near zero. We do this by induction on $d$, the degree of polynomial $P$.
For $\textrm{deg}(P)=1$, the claim is equivalent to the well known fact 
that the Fourier transform of a Calderon-Zygmund kernel is bounded.
For the inductive step put $P(t)=t^d+Q(t)$ and we have:
$$
\int_{-1}^{1}e^{iP(t)}k(t)dt=\int_{-1}^{1}[e^{it^d}-1][k(t)e^{iQ(t)}]dt+\int_{-1}^{1}e^{iQ(t)}k(t)dt$$

The second integral satisfies the condition by induction hypothesis,
and the first integral can be written in this form:
$$
\int_{-1}^{1}t^d\psi(t)k_{\epsilon,N}(t)dt$$ 
where $\psi$ is a smooth bounded function through the interval.
So the integrand is bounded, and integrated over a bounded interval it will produce a 
number that is bounded above, uniformly in $\epsilon$ and $N$. 
The integral near infinity is bounded, as can be seen by applying 
Van der Corput lemma of order $d$. $\square$


From this theorem it immediately follows that,
$m_{k_i}(\xi,\eta)$ is a uniformly bounded family of functions on the $(\xi,\eta)$ plane.
Just plug in $P(t)=2{\pi}(\xi t+ \eta t^2)$ in the previous theorem.


\medskip
\begin{lemma}
The limit
$m_{k}(\xi,\eta)
=\lim_{\epsilon \to 0, N \to \infty}m_{\epsilon, N}(\xi,\eta)$ 
exists and is a $C^{\infty}$ function for $\eta \neq 0$.
\end{lemma}


\noindent
{\bf Proof.}
This is rather standard but for the sake of completeness here is a proof.
We introduce ``parabolic'' polar coordinates $a,\delta$:
$(\delta,\delta^2)\otimes(a,1)=(\xi,\eta)$.
So for the upper half-plane, $\eta > 0$, we have:
$
\delta a=\xi, \delta^2=\eta, \qquad 
\textrm{or equivalently} \qquad a=\frac{\xi}{\eta^{1/2}}, \delta=\eta^{1/2}$
If $\eta \neq 0$, then $\delta \neq 0$ too, and it is obvious that 
the expression is smooth in $\delta$ on the region $\delta \neq 0$.
So there only remains the question of smoothness in $a$.
Notice that 
$m_{k}(\xi,-\eta)=m_{\widetilde{k}}(\xi,\eta)$ where $\widetilde{k}(x):=\bar{k(-x)}$.
So from now on we restrict our attention to $\eta > 0$ only.
We start with $k=k_{\epsilon,N}$. After a suitable change of variables we have:
$$
m_k(\xi,\eta)=\integral e^{-i(at+t^2)}k^{\delta}(t)dt
$$ 
Introducing smooth cut-off functions $\phi_1,\phi_2,\phi_3$, 
supported near $-\infty, 0, \infty$ respectively, we 
can write $m_k(\xi,\eta)$ as this sum:
$$
\int_{-1}^{1} e^{-i(at+t^2)}k^{\delta}(t)\phi_2(t)dt
+\int_{-\infty}^{-1} e^{-i(at+t^2)}k^{\delta}(t)\phi_1(t)dt+\int_{1}^{\infty} 
e^{-i(at+t^2)}k^{\delta}(t)\phi_3(t)dt
$$


First we treat the $\int_{1}^{\infty}$ part. The integral near $-\infty$ is similar.

\begin{displaymath}
\int_{1}^{\infty} e^{-i(at+t^2)}k^{\delta}(t)\phi(t)dt\\
=\int_{1}^{\infty}\frac{d}{dt} \{e^{-i(at+t^2)}\}  \frac{i}{a+2t} k^{\delta}(t)\phi(t) dt
=\int_{1}^{\infty}e^{-i(at+t^2)} \frac{d}{dt} \{ \frac{i}{a+2t} k^{\delta}(t)\phi(t)\} dt
\end{displaymath}

Continue integrating by parts and if we do it enough times we see that we have 
an absolutely convergent integral that decays in $a$ faster than any desired negative power of $a$,
so the expression is infinitely differentiable in $a$. Also it is clear from the definition that each of 
these integrals is infinitely differentiable in $\delta$ for $\delta \neq 0$.
The integral near zero is the Fourier transform of a compactly supported distribution so it is smooth.

Notice that all estimates are independent of $\epsilon$ and $N$, thus we can pass to the limit 
$\epsilon \to 0$, $N \to \infty$,  getting the result for $k$.
$\square$


\bigskip

We will now state a lemma, fully describing the asymptotics of the Fourier transform of $e^{ix^2}k(x)$, 
for an arbitary CZ kernel $k$. This lemma was generously provided to us by Elias Stein \cite{personal}.


\begin{lemma}
If $k$ is a Calderon-Zygmund kernel on $\R$, define
\begin{equation}
I_{k}(x):= \integral e^{-2{\pi}ixt}e^{-{\pi}it^2} k(t) dt
\end{equation}
then, $I_{k}(x)$ is a $C^{\infty}$ function of $x$
and it can be written as:
$I_{k}(x)=A(x)+e^{-i{\pi}x^2}B(x)$, 
where, for $x \to \infty$, we have: 
$A(x) \sim \widehat{k}(x) \qquad \qquad \textrm{and} \qquad \qquad B(x) \sim k(x)$

More precisely, for $ x \to \infty$, we have the following full asymptotic expansions:
$$
A(x) \sim \sum_{j=0}^{\infty}  {c_j} \partial_{x}^{2j} \widehat{k}(x) \qquad , \qquad
B(x) \sim \sum_{j=0}^{\infty} {c'}_j \partial_{x}^{2j} k(x)$$
\end{lemma}


{\bf Proof.}
Start with a partition of unity, separating the point zero, and infinity. Let $\eta$
be a $C^{\infty}$ function, that is identically equal to $1$ on $|x| \leq \frac{1}{2}$, 
and identically equal to zero, on $|x| \geq 1$.
Put $\eta_{0}:=\eta \quad \eta_{\infty}:=1-\eta_{0}$.
Decomposing $k$ accordingly:
$k=k_{0}+k_{\infty} \qquad \qquad k_{0}:=\eta_{0}k \quad , \quad k_{\infty}:=\eta_{\infty}k$


Observe that, on $|x| \geq C$, $\widehat{k_{0}}(x)$ is $C^{\infty}$, and 
$\widehat{k_{\infty}}(x)$ is a Schwartz function. 

That is because
$k_0$ is compactly supported so its Fourier transform is infinitely differentiable. 
Also, we have:
$
\widehat{k_{\infty}}(x)= \integral e^{-2{\pi}ixt}\eta_{\infty}(t)k(t) dt
$
But since $\eta_{\infty}$ is supported a distance away from zero, the integral is rapidly decreasing in $x$,
as can be seen by repeated integration by parts.


It follows from this observation that on the region $|x| \geq 1$ 
we have: $\widehat{k_0}(x)=\widehat{k}(x)+\phi$,
for some Schwartz function $\phi$; 
and similarly for $\widehat{k_{\infty}}(x)$ for $x$ small,  
$\widehat{k_{\infty}}(x)-\widehat{k}(x)$
is a $C^{\infty}$ function on $|x| \leq 1$.


First we observe that $I(x)$ is indeed a smooth function of $x$. 
Using the decomposition we just gave, put
$I_{k}(x)= I_{k_{\infty}}(x)+I_{k_{0}}(x)$. 
But $I_{k_{0}}$ is $C^{\infty}$ since it is the Fourier transform of a 
compactly supported distribution; 
just observe that on the ball of radius one, which is where $k_0$ is supported,
all $t^n e^{-i{\pi}t^2} k(t)$ are uniformly bounded.


For $I_{k_{\infty}}$ we have:
$$
I_{k_{\infty}}(x)= \int e^{-2{\pi}ixt} e^{-i{\pi}t^2}k(t)\eta_{\infty}(t)dt$$
which is a convergent integral. The argument goes as follows:

We have: $e^{-i{\pi}t^2}=\frac{i}{2{\pi}t}\frac{d}{dt}e^{-i{\pi}t^2}$,
so that
$$
I_{k_{\infty}}(x)=c \int \frac{1}{t}e^{-2{\pi}ixt}\{\frac{d}{dt}e^{-i{\pi}t^2}\}k(t)\eta_{\infty}(t)dt
$$


Integrating by parts we get a sum of four integrals, where the differentiation in $t$ falls on each of
the four terms $\frac{1}{t}$, $e^{-2{\pi}ixt}$, $k(t)$, and $\eta_{\infty}(t)$ respectively.
If the $t$ differentation falls on $e^{-2{\pi}ixt}$, we get an integral which is similar to 
the defining integral for $I_{k_{\infty}}$, but with a better decay. 
Repeating this process enough number of times, we will eventually get an absolutely convergent integral.

If the $t$-differentiation falls on $k(t)$ or $\frac{1}{t}$, we get an absolutely convergent integral, 
and if it falls on
$\eta_{\infty}$, the integral is the Fourier transform of a compactly supported distribution so 
it is a $C^{\infty}$ function.


Differentiating with respect to $x$ we have:
$$
\partial_{x} I_{k_{\infty}}(x)= c \int te^{-2{\pi}ixt} e^{-i{\pi}t^2}k(t)\eta_{\infty}(t)dt
=c \int e^{-2{\pi}ixt} \frac{d}{dt}e^{-i{\pi}t^2}k(t)\eta_{\infty}(t)dt
$$

and a carbon copy of the above argument gives the desired result on convergence of this integral.


Now there only remains the study of the asymptotic behaviour of this integral.


In the integral representing $I_k$, separate the contribution of $k_0$ and $k_{\infty}$:
$$
\integral e^{-2{\pi}ixt} e^{-i{\pi}t^2} k_{0}(t) dt + 
\integral e^{-2{\pi}ixt} e^{-i{\pi}t^2} k_{\infty}(t)dt$$

The first integral is the $A(x)$ in the statement of the lemma. 
Expand the term $e^{-i{\pi}t^2}$ in Taylor expansion around zero:
$
e^{-i{\pi}t^2}=\sum_{n=0}^{\infty} c_n t^{2n}$
and insert it in the integral. 
The desired expansion follows from the observation that 
$$
\int e^{-2{\pi}ixt} t^n k_{0}(t) dt = {c_n}\partial_{x}^{n} k(x)$$


The second term is equal to 
$$
e^{-i{\pi}x^2} \int e^{i{\pi}(x-t)^2} k_{\infty}(t) dt
$$
So we just need to show that this new integral has the asymptotic expansion we suggested for $B$.
\begin{equation}
B(x)=\int e^{i{\pi}(x-t)^2} k_{\infty}(t) dt = \int e^{i{\pi}t^2} k_{\infty}(x-t)dt
=c\int e^{-i{\pi}t^2} \widehat{k_{\infty}}(t) e^{-2{\pi}ixt} dt
\end{equation}

The last equality is true by unitariness of the Fourier transform and the fact that 
the Fourier transform of $k_{\infty}(x-t)$, is $\widehat{k_{\infty}}(t) e^{-2{\pi}ixt}$; 
and the Fourier transform of
$e^{-i{\pi}t^2}$ is $(i)^{-\frac{1}{2}}e^{i{\pi}x^2}$


In the last integral in equation 3.2, separate the contribution of zero and infinity, 
by bringing in the cut-off functions $\eta_0$ and $\eta_{\infty}$ again:
$$
\int e^{-i{\pi}t^2} \widehat{k_{\infty}}(t) e^{-2{\pi}ixt} dt = 
\int e^{-i{\pi}t^2} e^{-2{\pi}ixt} \eta_{\infty}(t) \widehat{k_{\infty}}(t) dt + I_{1}(x)$$

where $I_1$ is a similar integral where we have replaced $\eta_{\infty}$ by $\eta_0$.


But the first integral is rapidly decreasing since away from zero, you can integrate it by parts, 
as many number of times as desired, bringing down an order of decay in $x$ with each integration. 
So we only need to look at the integral $I_1$.
$$
I_{1}(x)= \int e^{-i{\pi}t^2} e^{-2{\pi}ixt} \eta_{0}(t) \widehat{k_{\infty}}(t) dt
$$
But as we observed in the observation above, near $x=0$, 
the term $\widehat{k_{\infty}}(x)$ is the sum of $\widehat{k}(x)$
and a $C^{\infty}$ function $\phi$. 
The remaining integral 
$$
\int e^{-i{\pi}t^2} e^{-2{\pi}ixt} \eta_{0}(t) \phi(t) dt$$ 
is easily seen to be of rapid decrease. 
The other integral is exactly the integral we had for $A(x)$, just 
with $\widehat{k}$ replacing $k$, so the same trick (Taylor expansion around zero) 
gives the similar result. 
Thus the proof of the lemma is now complete. $\square$


\begin{corollary}
If $\{ k_{\alpha} \}_{\alpha}$ is a uniformly bounded family of Calderon-Zygmund kernels,
and you define:
$$
I_{\alpha}(x):=\integral e^{-2{\pi}ixt} e^{-i{\pi}t^2} k_{\alpha}(t) dt$$

then each $I_{\alpha}$ is a smooth function of $x$, and
$$
I_{\alpha}(x)=A_{\alpha}(x)+e^{-i{\pi}x^2}B_{\alpha}(x)$$

where $A_{\alpha}$ and $B_{\alpha}$ satisfy the following differential inequalities
uniformly in $\alpha$:
$$
|\partial_{j} A_{\alpha}(x)| \leq {C_j}  |x|^{-j}   \qquad 
|\partial_{j} B_{\alpha}(x)| \leq {C'_j} |x|^{-1-j}  $$
\end{corollary}

This is clear from the proof of Lemma 3.2.


\vspace{6mm}

\section{Proof of theorem 1.1}

In the following we will give a proof of the theorem for the case ${c_0}=1$.
The proof for arbitary constant ${c_0}$ is no different.

%
%


First observe that, by corollary 2.4.4 of Nagel, Ricci, Stein \cite{flag}.
for every flat kernel $M$, there is a uniformly bounded collection of compactly supported, 
$C^{\infty}$ functions $\{ \phi_{m,n} \}_{m,n \in \mathbb{Z},n \geq 0}$,
whose supports are bounded uniformly and
which satisfy the estimates for flag kernels uniformly,
and
each $\phi_{m,n}$ satisfies the cancellation conditions:
$$
\int \phi_{m,n}(x,y)dx=0 \qquad \int \phi_{m,n}(x,y)dy=0$$

such that $\sum_{m,n} \phi_{m,n}^{(m,m+n)}$ converges to $M$ in the sense of distributions.
$\phi^{(m,n)}$ is a dyadic re-scaling of $\phi$ defined as:
$
\phi^{(m,n)}(x,y):= 2^{-m-n}\phi(2^{-m}x,2^{-n}y) \qquad \textrm{dyadic re-scaling}
$


This gives a sequence $\{f_i\}$ of $C^{\infty}_{0}$ functions which converges, 
in the sense of distributions, to $M$, and which satisfy uniformly the estimates for flag kernels.

As a result, the $C^{\infty}_{0}$ functions $\tilde{f}_i(x,y)=K_i(x,y):=f_i(y-cx^2)$ converge 
in the sense of distributions to the curved kernel $K$, and satisfy uniformly the estimates for $K$.
But we know that if a uniformly bounded sequence of distributions $\{ f_{i} \}$ converges to a distribution $F$,
then $\{ \widehat{f}_{i} \}$ converges to $\widehat{F}$ in the sense of distributions. 
Thus the functions $\widehat{K_i}$
converge in the sense of distributions to $\widehat{K}$ and by what will be proved below 
(for $K_i$ in the place of $K$) the functions $\widehat{K_i}$ will then be uniformly bounded. 
Hence $\widehat{K}$ will be a bounded function. 
Our proof then proceeds with $K$ standing for $K_i$, and $M$ standing for $f_i$.



If $M_{\eta}$ is the Fourier transform in the second variable:
$$
M_{\eta}(x):= \int e^{-2{\pi}iy\eta} M(x,y) dy$$
we observe that
\begin{equation}
m(\xi,\eta)=\int e^{-2{\pi}i(x\xi+x^2\eta)} M_{\eta}(x)dx
\end{equation}



It is easy to check that $\{M_{\eta}\}_\eta$ is a uniformly bounded family of CZ kernels. 
Let us for example check the cancellation condition and see if they are 
uniformly satisfied:
$$
\int M_{\eta}(x)\phi(Rx)dx = \int\!\!\int e^{-2{\pi}iy\eta} M(x,y)\phi(Rx) dxdy
$$
but by definition the family $\int M(x,y) \phi(Rx) dx$ is a uniformly bounded family of
CZ kernels, so the bounds on their multipliers are uniform and we are done.


We sometimes wish to write our coordinates in polar form. 
Whenever we have a function, distribution, etc, in $\xi,\eta$ variables (e.g. $M(\xi,\eta)$)
we denote the same quantity in polar coordinates, by the same symbol, with a $\pi$ subscript:
$$
M_{\pi}(a,\delta):=M(a\delta,\delta^2)=M(\xi,\eta)
$$


Switching to this system of polar coordinates, we are ready to finish the proof of 
thorem 1.1 which claims that each of $L_{1,\pi}$ and $L_{2,\pi}$ brings down one order of
decay in $\delta$ with each differentiation, and one order of decay in $a$ for large $a$ 
( $a \to \infty$).
$$
|\partial^{\alpha}_{\delta} \partial^{\beta}_{a} L_{i,\pi} (a,\delta)| \leq C_{\alpha,\beta} |\delta|^{-\alpha} (1+a)^{-\beta}
$$
This is equivalent to the claim that each order of differentiation
in $a$ of each of the uniformly bounded families of CZ multipliers
$$
\delta^{\alpha} \partial_{\delta}^{\alpha}L_{1,\pi}(a,\delta) 
\hspace{2cm}
\delta^{\alpha} \partial_{\delta}^{\alpha}L_{2,\pi}(a,\delta) 
$$
brings down a decay of $|a|^{-1}$ for $a \to \infty$. This last claim is 
the content of Corollary 2.4 on the asymptotic behaviour
of a certain integral defined for a uniformly bounded family of Calderon-Zygmund kernels.
So, now we only need to show that on the region where $|\eta|>c|\xi|^2>0$ the multiplier of the curved kernel 
is a Mikhlin multiplier and that will end the proof.


So we know the behaviour of $m(\xi,\eta)$ on the region where $|\xi| \geq c|\eta|^{1/2}$. There only remains 
the analysis of this Fourier transform in the region above the parabola: $|\eta| > c|\xi|^2 > 0$.


\begin{lemma}
For each $\alpha$, the family 
$\left\{ \eta^{\alpha} \partial^{\alpha}_{\eta} M_{\eta}(x) \right\}_{\eta}$
is a uniformly bounded family of CZ kernels. 
\end{lemma}

{\bf Proof.}
If $M(x,y)$ is a flag kernel with respect to flag $\mathcal{F}_1$, then it is easy to check that 
the partial derivative of $M$ with respect to the second variable $y$, 
is a new distribution $M'$ which is itself a flag kernel, with respect to the same flag,
i.e. for every $M$, there exists an $M'$ such that $y\partial_y M(x,y)= M'(x,y)$.
But then we have:
$$
\eta \partial_{\eta} M_{\eta}(x)= \int \eta e^{-iy\eta}yM(x,y) dy
=\int \partial_{y}\{ e^{-iy\eta}\} M(x,y) dy= \int e^{-iy\eta} \partial_{y}\{yM(x,y)\} dy$$


\bigskip

\begin{lemma}
On the region where $|\eta| > c|\xi|^2 > 0$, the Fourier transform of $K$ is
a Mikhlin multiplier for non-isotropic dilations:
$$
|\partial_{\xi}^{\alpha} \partial_{\eta}^{\beta} m(\xi,\eta)| 
\leq C_{\alpha,\beta} (|\xi|+|\eta|^{1/2})^{-\alpha-2\beta}$$

\end{lemma}


{\bf Proof.} 
We need to study the behaviour of the integral representing $m(\xi,\eta)$
under differentiation with respect to $\xi$ and $\eta$, and check whether they satisfy
the decay properties needed. 
We do this, by checking separately the decay condition 
for the first derivatives with respect to $\xi$ and $\eta$, 
and the similar conclusion for higher derivatives will follow. 
Namely, we will check whether
$$
|\partial_{\xi}m(\xi,\eta)| \leq C(|\xi|+|\eta|^{1/2})^{-1} \quad \textrm{and} \quad
|\partial_{\eta}m(\xi,\eta)| \leq C(|\xi|+|\eta|^{1/2})^{-2}
$$

We have:
$$
\partial_{\xi} m(\xi,\eta)=\int {x} M_{\eta}(x) e^{-i(x\xi+x^2\eta)}dx$$

Break-up this integral into two parts: $|x| < A |\eta|^{1/2}$ and $|x| \geq A |\eta|^{1/2}$. ($A$ to be specified later)\\
For $|x| < A |\eta|^{1/2}$ use the fact that no matter what $A$ is 
$\{M_{\eta}\}$ and 
$\{\eta^{-\frac{1}{2}}xM_{\eta}(\eta^{-\frac{1}{2}}x)\}_{\eta}$
are uniformly bounded families of Calderon-Zygmund kernels (lemma 3.2). 
So, we get
$$
|\int \eta^{-\frac{1}{2}}x M_{\eta}(\eta^{-\frac{1}{2}}x) e^{-ix\frac{\xi}{\eta^{1/2}}}
e^{-ix^2}dx| \leq C
$$
for some constant C, not depending on $\xi$ or $\eta$.

Put $x':=\eta^{-\frac{1}{2}}x$ and we get:
$$
|\int x'M_{\eta}(x') e^{-ix'\xi} e^{-i{x'}^{2}\eta} dx'| \leq C |\eta|^{-\frac{1}{2}}
$$
which is the same as 

\begin{equation}
|\partial_{\xi} m(\xi,\eta)| \leq C |\eta|^{-\frac{1}{2}}
\end{equation}


But we are in the region where $\eta$ is comparatively big ($|\eta|^{1/2} \geq c |\xi|$), 
so the dominating factor in $(|\xi|+|\eta|^{1/2})$
is $|\eta|^{1/2}$, that is: 
$c_1|\eta|^{1/2} \leq |\xi|+|\eta|^{1/2} \leq c_2 |\eta|^{1/2}$ for constants $c_1,c_2$.
This means that inequality (5) is equivalent to the inequality:
$$
|\partial_{\xi} m(\xi,\eta)| \leq C (|\xi|+|\eta|^{\frac{1}{2}})^{-1}$$
which is exactly what we claimed.


to treat $|x| \geq A|\eta|^{1/2}$, 
pick $A$ such that the stationary point of the phase ($x_0=\frac{2\xi}{\eta}$)
falls outside the support of integral $|x| \geq A|\eta|^{1/2}$,
therefore reducing it to a nonstationary phase integral which is rapidly decreasing. 
This takes care of differentiation in $\xi$.


For treatment of the case of differentiation with respect to $\eta$, 
notice that:
\begin{equation}
\partial_{\eta} m(\xi,\eta)=
\int \{ \partial_{\eta}M_{\eta}(x)+M_{\eta}(x)x^2 \} e^{-i(x\xi+x^2\eta)} dx
\end{equation}


But the contribution from the first term of the integrand in equation (3.3) is taken care of
by corollary 2.4, showing that it brings down one order of decay in $\eta$,
which is exactly what we need.


The second integral in equation (3.3),
also brings down an order of decay in $\eta$. 
The argument 
goes exactly like the case for differentiation in $\xi$ that resulted
a decay of $|\eta|^{-\frac{1}{2}}$ in equation (3.2).
The proof of Theorem 1.1 is now complete. $\square$


Using theorem 1.1, we can also deduce a boundedness result for the class of curved flag kernels under study.
This was independently proved by Secco \cite{secco}. After our work was completed in 2004,
we found out about the paper of Secco, which is a direct proof of this corollary.


\begin{corollary}
If $M$ is a flat flag kernel with respect to flag $\mathcal{F}_1$, and $K$ is the curved version,
and $m$ the Fourier transform of $K$, then $m$ is an $L^p$ multiplier for $1<p<\infty$.
\end{corollary}


\noindent
{\bf Proof.}
Introducing a smooth cut-off function, break up $m$ into two parts.\\
Suppose $\Psi(x)$ is a compactly supported $C^{\infty}$ function 
which is identically equal to 1 on $|x| \leq 1$, 
and zero on $|x| >2$. Put:
$$
m_0(\xi,\eta):= m(\xi,\eta) \Psi(\frac{\eta}{\xi^2}) \qquad m_1:=m-m_0
$$
$m_0$ being the restriction on the region under parabola of a fixed slope: 
$|\xi| \geq c|\eta|^{1/2}$, 
and $m_1$ restriction on the region above the parabola: $|\xi| < c|\eta|^{1/2}$.
$m_1$ is a Mikhlin multiplier and so an $L^p$ multiplier, as can be seen from theorems of 
Fabes, Riviere, \cite{studia-66}, and Folland, Stein. \cite{folland}.
$m_0$ is an $L^p$ multiplier because it is the sum of a flat flag multiplier which we know is
an $L^p$ mutplier \cite{flag}, and a second term 
which is itself a product of two terms.
$m_0(\xi,\eta)=L_1(\xi,\eta)+e^{ic'\frac{\xi^2}{\eta}} \frac{\eta^{1/2}}{\xi} L_2(\xi,\eta)$
One term is again a flat flag mutiplier so we use 
the theorem of Nagel, Ricci, and Stein again. 
The other term is basically the multiplier of 
ordinary Hilbert transform on the parabola, so it is bounded on $L^p$. $\square$


\bigskip

By analogy to the special case studied before, 
we can predict the asymptotic behaviour of
$L_1$ and $L_{2}$:
$$ 
\textrm{when }\frac{\xi}{\eta^{1/2}} \to \infty:
\hspace{2cm}
L_{1}(\xi,\eta) \sim \widehat{M}(\xi,\eta),
\qquad
L_{2}(\xi,\eta) \sim \frac{\xi}{\eta}M_{\eta}(\frac{\xi}{\eta})
$$

In fact, using lemma 2.4 we can establish a more precise relation between $L_1$ and $L_2$.
This is the content of our remark in the first section.
Howevere, we note that the relation can only be stated in terms of asymptotics 
of these multipliers, since the exact form of the multiplier depends on many factors such as 
the specific cut off functions we use at different stages and hence out of control.

\bigskip
Notice that the term $\widehat{M}(\xi,\eta)$ is also a flat flag multiplier, but 
with respect to the flag $\mathcal{F}_2$ instead. So on the region specified, 
it is highly regular: it is a Mikhlin multiplier
with respect to the specified non-isotropic dilations.
Also notice that $\frac{\eta^{1/2}}{\xi}L_2(\xi,\eta)$ is a flag multiplier with respect to 
flag $\mathcal{F}_1$ on $|\xi| \geq c|\eta|^{1/2}$.


\bigskip

The constant $c'$ is related to the constant $c_0$ in choosing the parabola that carries 
the singularity, 
by ${c_0}c'=\frac{\pi}{2}$. 
This can be seen from the following heuristic argument. The Fourier transform of $K$ is such an integral:
$$
\widehat{K}(\xi,\eta)=\int\!\!\int e^{-2{\pi}i(x\xi+y\eta)} M(x,y-{c_0}x^2) dxdy
$$
$$
=\int\!\!\int e^{-2{\pi}i(x\xi+{c_0}x^2\eta)}e^{-2{\pi}iy\eta}M(x,y)dxdy
=\int e^{-2{\pi}i(x\xi+{c_0}x^2\eta)}M_{\eta}(x)dx
$$
for $\{ M_{\eta} \}$ a uniformly bounded family of CZ kernels.


The main contribution in this last integral is from the stationary phase point
at $x_0=\frac{-\xi}{2{c_0}\eta}$. We know that the main term in a stationary integral of the form
$\int e^{i\Phi(x)} \Psi(x)dx$ has a phase term $e^{i\Phi(x_0)}$ 
where $x_0$ is the stationary point.
So for our integral the phase term is $e^{-2{\pi}i\frac{-\xi^2}{4{c_0}\eta}}$, 
which is the same as $e^{ic'\frac{\xi^2}{\eta}}$.


\section{proof of theorem 1.2}

The machinery we need has already been developed for the proof of theorem 1.1, and
the proof of the reverse direction is basically a bootstrap argument.
First we observe that

\begin{lemma}
If $M$ is a flat flag kernel with respect to flag $\mathcal{F}_2$, then 
for any constant $c$, the new kernel defined as $M(x,y-cx^2)$ is also 
a flat flag kernel with respect to the same flag.
\end{lemma}

This is a rather trivial lemma, remembering the alternative description of the kernels
and multipliers in parabolic polar coordinates.
When you shift to parabolic polar coordinates, the change of variables described
in the satement of the lemma is the same as a ``shift'' in $a$ (parabolic angle)
so none of the conditions change.

{\bf Proof of theorem 1.2.}
Writing down the inverse Fourier transform of the given distribution, it will be of the form:
\begin{equation}
K(x,y):= 
{\tt Inverse  F.T. }
\{\Phi(\xi,\eta)e^{ic'\frac{\xi^2}{\eta}}\frac{\eta^{1/2}}{\xi}\ell(\xi,\eta)\}(x,y)
=
\int\!\!\int e^{2{\pi}ix\xi} e^{2{\pi}iy\eta} \Phi(\xi,\eta) e^{ic'\frac{\xi^2}{\eta}} \frac{\eta^{1/2}}{\xi} \ell(\xi,\eta) d\xi d\eta
\end{equation} 
Which we have to show is in this form: $M_2(x,y)+M_1(x,y-{c_0}x^2)$ where $M_1$ is a flat flag kernel with respect to $\mathcal{F}_1$
and $M_2$ is a flat flag kernel with respect to $\mathcal{F}_2$ and $c_0$ the dual constant to $c'$.
Or, equivalanatly, that $K(x,y+{c_0}x^2)$ is of the form $M_1(x,y)+M_2(x,y+{c_0}x^2)$
But according to lemma 4.1 $M_2(x,y+{c_0}x^2)$ is just $M'(x,y)$ for some other flat flag kernel $M'$
with respect to the same flag, and to show that $K(x,y+{c_0}x^2)$ is in form $M_1(x,y)+M'(x,y)$
where $M_1$ and $M'$ are flat flag kernels with respect to $\mathcal{F}_1$ and $\mathcal{F}_2$
respectively, is the same as showing that it equals a product kernel.

Performing integration in (4.1) with respect to $\xi$ first, the distribution can be written as the following double integral:
\begin{equation}
K(x,y)=
\int e^{2{\pi}iy\eta} \eta^{1/2} \{ \int_{|\xi|^2 \geq c |\eta|} e^{2{\pi}ix\xi} \Phi(\xi,\eta) e^{ic'\frac{\xi^2}{\eta}} \frac{1}{\xi} \ell(\xi,\eta)d\xi \} d\eta
\end{equation}

The main contribution is from the stationary point of the phase at $\xi_0=cx\eta$, 
so by absorbing the error in the $\Phi$ term, replacing it with another cut off function $\tilde{\Phi}$,
satisfying the same decay estimates, the whole expression reduces to
$$
\int_{|x^2 \eta| \geq c} e^{2{\pi}iy\eta} \eta^{1/2} e^{i{c_0}x^2\eta} \left(\frac{2\pi}{-i2c'{\eta}^{-1}}\right)^{1/2} 
\tilde{\Phi}(x,\eta) \frac{1}{cx\eta} \ell(cx\eta,\eta) d\eta
$$

but according to lemma 4.1, we only need to show that this expression equals a product kernel after the change of variables 
$y \mapsto y-{c_0}x^2$.
This is the content of next lemma:

\begin{lemma}
$$
P(x,y):= K(x,y+{c_0}x^2)=\int_{|x^2 \eta| \geq c} c  e^{2{\pi}iy\eta} \tilde{\Phi}(x,\eta) \frac{1}{x} \ell(cx\eta,\eta) d\eta
$$
is a product kernel in $\mathbb{R}^2$.
\end{lemma}

We prove the claim by taking the Fourier transform of this expression and show that it is a bounded function
and has the decay properties of a product multiplier and that would conclude the proof.
$$
\widehat{P}(\xi,\eta)=
\int_{|x^2{\eta}| \geq c} e^{-2{\pi}ix\xi}\frac{1}{x}\ell(x\eta,\eta)\tilde{\Phi}(x,\eta) dx =
\int_{|x^2{\eta}| \geq c} e^{-2{\pi}ix\xi}\frac{1}{x} \tilde{\ell}(x\eta,\eta) dx 
$$
for another flat flag multiplier of the same kind $\tilde{\ell}(\xi,\eta)$.
It is clear that away from the coordinate axes $\eta=0$ and $\xi=0$
this expression is a smooth function.

But $\tilde{\ell}$ is a flag multiplier and hence also a product multiplier, so 
$$
\partial_{\eta}\tilde{\ell}(x\eta,\eta)=x\partial_1 \ell(x\eta,\eta) + \partial_2 \ell(x\eta,\eta)
=x(x\eta)^{-1}{\ell_1}(x\eta,\eta)+\eta^{-1}\ell_2(x\eta,\eta)$$
for $\ell_1$ and $\ell_2$ of the same type.
So each differentiation in $\eta$ brings down one order of decay $|\eta|^{-1}$ and leaves 
an integral of the same kind. So there only remains treatment of differentiation in $\xi$.

Notice that the integral is in form of the Fourier transform of a multiple of two functions,
namely $\frac{1}{x}$ (restricted to $|x^2\eta| \geq c$) and $\tilde{\ell}(x\eta,\eta)$.
So, the integral will be the convolution of their Fourier transforms. Also notice that the Fourier transform 
of the first function $\frac{1}{x}$ is uniformly bounded, so we only need to show that
$$
|\partial_{\xi}^{\alpha} \int_{|x^2\eta| \geq c} e^{-2{\pi}ix\xi} \tilde{\ell}(x\eta,\eta) dx| 
\leq C_{\alpha} |\xi|^{-\alpha}
$$

Break up the integral into two parts:
$$
\left( \int_{|x\xi| < 1} + \int_{|x\xi| \geq 1} \right) e^{-2{\pi}ix\xi} \tilde{\ell}(x\eta,\eta) dx 
$$

For the fist integral we want to show that 
$$|\int_{|x\xi| < 1} x\xi e^{-2{\pi}ix\xi} \tilde{\ell}(x\eta,\eta) dx| \leq C $$
which is trivial since $\tilde{\ell}$ is bounded. 
For the second integral we want to show that
$$|\int_{|x\xi| < 1} \xi e^{-2{\pi}ix\xi} \tilde{\ell}(x\eta,\eta) dx| \leq C |x|^{-1}$$
which is clear after an integration by parts. (This is basically a nonstationary phase estimate decay)
Higher order derivatives are treated the same way.
Thus the proof of theorem 1.2 is now complete. $\square$

\bigskip

A remark about asymptotic behaviour of $M_1$ and $M_2$ is in order.
The inverse Fourier transform maps $\Phi(\xi,\eta)e^{ic'\frac{\xi^2}{\eta}}\eta^{1/2}\xi^{-1}\ell(\xi,\eta)$
to $M_1(x,y)+M_2(x,y-c_0x^2)$. Taking the Fourier transform back again, we should end up with the original function.
Yet, using theorem 1.1 there is an alternative description for this expression as a sum of three terms:
$$
\widehat{M_1}(\xi,\eta)+L_1(\xi,\eta)+\Psi(\xi,\eta)e^{ic'\frac{\xi^2}{\eta}}\eta^{1/2}\xi^{-1}\ell(\xi,\eta)$$

In the limit where $\frac{\xi^2}{\eta} \to 0$, $\Psi$ and $\Phi$ are identically equal to $1$
and hence $\widehat{M_1}(\xi,\eta) \sim L_1(\xi,\eta)$.
But by the remark from theorem 1.1 we know that $L_1(\xi,\eta) \sim \widehat{M_2(}\xi,\eta)$
which concludes the remark on theorem 1.2.


\bigskip

{\bf Aacknowledgements}

The major part of this work is from the author's dissertation
written under the guidance of Professor Elias Stein,
whom we gratefully thank for many helpful discussions 
as well as invaluable insight.
The author was partially supported by NSERC grant number 22R80520.



\medskip
\noindent
\sc{ Department of Mathematics and Statistics, University of Helsinki, Finland}

\noindent
\it {Hadi Jor\'{a}ti, hadi.jorati@helsinki.fi}



\newpage

{\small



\begin{thebibliography}{99}

\bibitem{annals} M. Christ, A. Nagel, E. M. Stein, S. Wainger,, 
{\it Singular and Maximal Radon transforms: Analysis and geometry}, 
Annals of Math. {\bf 150} (1999), 489-577.


\bibitem{studia-66} E. B. Fabes, N. M. Riviere,
{\it Singular integrals with mixed homogeneity},
Studia Mathematica, T. XVII (1966), 19-38






\bibitem{flag} A. Nagel, F. Ricci, E. M. Stein, 
{\it Singular integrals with flag kernels and analysis on quadratic CR manifolds}, 
J. Functional Analysis {\bf 181} (2001), 29-118.



\bibitem{pre-2} A. Nagel, E. M. Stein, 
{\it On the product theory of singular integrals}, 
Revista Math. Iberoamericana {\bf 20} (2004), 531-561.



\bibitem{secco} S. Secco,
{\it Adapting product kernels to curves in the plane},
Math. Zeit. {\bf 248} (2004), 459-476.



\bibitem{personal} E.M. Stein
personal communication.


\bibitem{bigstein} E.M. Stein,
{\it Harmonic Analysis},
Princeton University Press, 1993


\bibitem{kahane} E. M. Stein, 
{\it Oscillatory integrals related to Radon-like transforms}, 
The Journal of Fourier Analysis and Applications (Kahane special issue) (1995), 535-551.


\bibitem{studia-70} E. M. Stein, S. Wainger,
{\it The estimation of an integral arising in multiplier transformations},
Studia Mathematica, T. XXXV. (1970), 101-104


\bibitem{bull} E. M. Stein, S. Wainger, 
{\it Problems in Harmonic Analysis related to curvature}, 
Bull. AMS {\bf 84} (1978), 1239-1295.



\end{thebibliography}
}
\end{document}